\documentclass[12pt]{article}

\usepackage{amsmath}
\usepackage{amsthm}
\usepackage{amsfonts}
\usepackage{mathrsfs}
 \usepackage{stmaryrd}
\usepackage{setspace}
\usepackage{fullpage}
\usepackage{amssymb}
\usepackage{breqn}
\usepackage{enumitem}
\usepackage{bbold}
\usepackage{authblk}
\usepackage{comment}
\usepackage{hyperref}
\usepackage{pgf,tikz}
\usepackage{graphicx}
\usepackage{subcaption}

\newtheorem{thm}{Theorem}[section]

\newtheorem{lemma}[thm]{Lemma}
\newtheorem{conjecture}[thm]{Conjecture}
\newtheorem{proposition}[thm]{Proposition}
\newtheorem{prop}[thm]{Proposition}

\newtheorem{corollary}[thm]{Corollary}

\newtheorem{clm}[thm]{Claim}

\newcommand\ex{\ensuremath{\mathrm{ex}}}

\newcommand\cF{{\mathcal F}}
\newcommand\cG{{\mathcal G}}
\newcommand\cH{{\mathcal H}}

\newtheorem*{thm*}{Theorem}
\newtheorem*{prop*}{Proposition}
\newcommand{\ignore}[1]{}

\title{Multicolor hypergraph Tur\'an problems}
\author{
Xin Cheng\thanks{\small School of Mathematics and Statistics, Northwestern Polytechnical University and Xi'an-Budapest Joint Research Center for Combinatorics, Xi'an 710129, Shaanxi, P.R. China. Email:
\small \texttt{xincheng@mail.nwpu.edu.cn}.}\,, \hspace{0.2em}
D\'{a}niel Gerbner\thanks{\small Alfr\'ed R\'enyi Institute of Mathematics and Department of Computer Science and Information Theory, Faculty of
Electrical Engineering and Informatics, Budapest University of
Technology and Economics, Műegyetem rkp. 3., H-1111 Budapest, Hungary. Email:
\small \texttt{gerbner.daniel@renyi.hu}.}\,, \hspace{0.2em}
Yahui Zhang\thanks{\small
School of Mathematics, Faculty of Basic Sciences,
University of Shanghai for Science and Technology,
Shanghai 200093, China;
Alfr\'ed R\'enyi Institute of Mathematics,
Budapest, Hungary. Email:
\small \texttt{241610179@st.usst.edu.cn}.}}

\date{}

\begin{document}

\maketitle

\begin{abstract}
We study rainbow Tur\'an problems for collections of $k$ $r$-uniform hypergraphs on a common $n$-vertex set. For a family $\mathcal{F}$ of forbidden $r$-graphs, we establish a general upper bound for the maximum total number of hyperedges in a rainbow-$\mathcal{F}$-free collection. 
This yields an exact formula when $k$ is sufficiently large. 
We also prove that, for every $r$-partite $r$-graph $F$ with $\ell$ edges, the rainbow Tur\'an number equals $(\ell-1)\binom{n}{r}$ for all sufficiently large $n$. 
For the $3$-graph $F_5$, we determine the rainbow Tur\'an number for every fixed $k\geq 3$ and sufficiently large $n$. 
Finally, for expansions $F^{r+}$ of color-critical graphs $F$ with $\chi(F)=t+1>r$, we show that the rainbow Tur\'an number is $k|E(T^r(n,t))|$ whenever $k$ and $n$ are sufficiently large.

\end{abstract}



\section{Introduction}

One of the most fundamental results in extremal Combinatorics is Tur\'an's theorem \cite{Tu}, which determines the largest number of edges in an $n$-vertex $K_k$-free graph. Generalizations of this problem have attracted a lot of attention.

Given a family  $\mathcal{F}$ of $r$-uniform hypergraphs, we denote by $\ex^r(n,\cF)$ the largest number of hyperedges in an $n$-vertex $r$-uniform hypergraph that does not contain any member of $\mathcal{F}$ as a subhypergraph. If $\mathcal{F}=\{F\}$, we use the simpler notation $\ex^r(n,F)$ instead of $\ex^r(n,\{F\})$.

Hilton \cite{hilt} and, independently, Keevash, Saks, Sudakov and Verstra\"ete \cite{Keevash2004} introduced a multicolor variant. 
Given a collection of $k$ graphs $G_1,\dots,G_k$ on the same vertex set, a \textit{rainbow} copy of a graph $F$ is a copy of $F$ in which each edge is assigned to a graph in the collection containing that edge, such that distinct edges are assigned distinct graphs. 
Note that an edge can belong to multiple graphs $G_i$, and we have to pick one of them.

Here we consider the hypergraph version of this problem. We denote by $\ex^r_k(n,\mathcal{F})$ the largest sum of the number of hyperedges contained in a collection of $k$ $r$-uniform hypergraphs on the same $n$-element vertex set without a rainbow copy of any member of $\mathcal{F}$. 

To the best of our knowledge, this problem for $r>2$ has been studied only in some special cases. For two independent edges, this is called the cross-intersecting property and it is well-studied, especially in the case of two colors. More complicated forbidden intersection patterns (i.e., forbidding a family of hypergraphs each with two hyperedges) can be found in \cite{ks}, with sharper bounds in the case of two colors in \cite{ks2}. Larger matchings were studied in \cite{hilton} and also in \cite{bks} where some other hypergraphs (including non-uniform ones) were also studied. 
Zhao, Wang and Zhou \cite{zwz} proved some exact results when $\cF$ is the expansion of a graph clique and used them to deal with an uncolored hypergraph Tur\'an problem.

Let us describe two simple rainbow $\cF$-free constructions. Let $\ell$ be the smallest number of hyperedges of a member of $\cF$, then the collection of $\ell-1$ complete hypergraphs and $k-\ell+1$ empty hypergraphs is clearly rainbow $\cF$-free, giving the lower bound $\ex^r_k(n,\mathcal{F})\ge (\ell-1)\binom{n}{r}$. Let $G$ be an $\cF$-free $r$-graph with $\ex^r(n,\cF)$ hyperedges, and let $G_1=G_2=\dots=G_k=G$. This collection is also $\cF$-free, giving the lower bound $\ex^r_k(n,\mathcal{F})\ge k\ex^r(n,\cF)$.

In the graph case, Keevash, Saks, Sudakov and Verstra\"ete \cite{Keevash2004} proved for several graphs and integers $k$ that one of these two lower bounds holds with equality, and showed an example where this is not even asymptotically true. 
In particular, their Theorem 1.1 states that if $k\ge \binom{n}{2}-\ex^2(n,F)+\ell$, then $\ex^2_k(n,F)= k\ex^2(n,F)$. Our first result shows that this phenomenon extends to hypergraphs.

\begin{thm}\label{thm: sum of edges non-r-partite large k} Let $F$ be an $r$-graph and let $\ell = e(F)$.
        If $k\ge (\ell-1)\binom{n}{r}-\ell\ex^r(n,F)+\ell$, then
  $ex_k^r (n, {F}) = k ex^r(n, F)$.
\end{thm}

Note that in the case $r=2$, our result improves the threshold in Theorem 1.1 of \cite{Keevash2004} when $(\ell-2) \binom{n}{2} < (\ell -1)ex^2(n, F)$. The above result easily follows from the theorem below.

    \begin{thm}\label{thm: sum of edges non-r-partite}
        Let $r \geq 2$, $\ell \geq 2$ and $k \geq \ell$.
        Let $\mathcal{F}$ be a family of $r$-uniform hypergraphs such that $\ell$ is the smallest number of hyperedges in a member of $\cF$.
        Let $F \in \cF$ with $\ell = e(F)$.
        If $n \geq v(F)$, then we have

        \begin{equation*}
            \ex_k^r (n, \mathcal{F}) \leq (\ell - 1) \binom{n}{r} + (k - \ell) \ex^r(n, F).
        \end{equation*}
    \end{thm} 

    Another simple corollary is the following.

    \begin{corollary} For any $\cF$ such that $\ell$ is the smallest number of hyperedges in a member of $\cF$, we have
        $ex_\ell^r(n, \mathcal{F}) = (\ell - 1) \binom{n}{r}$ for $n \geq v(F)$ where $F \in \cF$ with $\ell = e(F)$.
    \end{corollary}

    Let us continue with a general lower bound. Let $\cF(i)$ denote the family of hypergraphs obtained by deleting $i$ hyperedges from a member of $\cF$.

\begin{proposition}
    We have $ \ex_k^r (n, \mathcal{F})\ge \max \left\{i \binom{n}{r} + \ex_{k-i}^r(n, \cF(i)): i\ge 0\right\} $.
\end{proposition}

Note that the proof of the above proposition is trivial: we take $k-i$ hypergraphs without a rainbow $\cF(i)$ and $i$ complete hypergraphs.
The resulting collection contains no rainbow member of $\cF$, since any such rainbow copy of some $F\in\cF$ would contain
a rainbow subhypergraph with $e(F)-i$ hyperedges in the first $k-i$ hypergraphs.
But this is a rainbow member of $\cF(i)$, a contradiction.
Also note that this construction contains the two simple constructions at $i=0$ and $i=\ell-1$. It also contains the example we mentioned earlier where the two simple constructions are not extremal (in the graph case). \cite{chak} contains further examples; those are also contained here. 

\begin{conjecture}
    We have $ \ex_k^r (n, \mathcal{F})= \max \left\{i \binom{n}{r} + \ex_{k-i}^r(n, \cF(i)): i\ge 0\right\} $ for sufficiently large $n$.
\end{conjecture}

Note that if $H_1$ is the complete hypergraph on $n$ vertices, then the rest of the collection must avoid a rainbow member of $\cF(1)$, thus has at most $\ex_{k-1}^r(n, \cF(1))$ hyperedges. If the conjecture holds for $\cF(1)$ and $k-1$, then this upper bound is equal to $\max \left\{(i-1) \binom{n}{r} + \ex_{k-i}^r(n, \cF(i)): i\ge 1\right\}$, thus the conjecture also holds for $\cF$ and $k$. The more problematic case is where $H_1$ is not complete. According to the conjecture, $H_1$ is an extremal $\cF$-free hypergraph in this case.

    Gerbner and Miao \cite{gemi} proved that for any bipartite graph $F$ with $\ell$ edges and positive integer $k \geq \ell$, for sufficiently large $n$ we have $ex_k^2(n, {F}) = (\ell - 1) \binom{n}{2}$. Our next theorem extends this to hypergraphs.

    \begin{thm}\label{thm: sum of edges r-partite}
        Let $r \geq 2$, $\ell \geq 2$ and $k \geq \ell$.
        Let $F$ be an $r$-uniform $r$-partite hypergraph with $\ell$ hyperedges.
        For fixed $k$, when $n$ is sufficiently large, we have
        \begin{equation*}
            \ex_k^r(n, F) = (\ell - 1) \binom{n}{r}.
        \end{equation*}
        The bound is achieved by taking $\ell-1$ copies of $K_n^r$ and $k-\ell+1$ copies of empty hypergraphs.
    \end{thm}

    Recall that several graph results in this area show that for every $k$, one of the two simplest constructions is extremal, for sufficiently large $n$. Here we give a hypergraph example of this phenomenon.

    Let $F_5$ denote the 3-graph with hyperedges $\{1,2,3\}, \{1,2,4\}, \{3,4,5\}$. 
    Let $T^r(n,t)$ denote the complete $t$-partite $r$-uniform hypergraph with each part of size $\lfloor n/t\rfloor$ or $\lceil n/t\rceil$. Frankl and F\"uredi \cite{ff} proved that $\ex^3(n,F_5)=|E(T^3(n,3))|$ for $n\ge 3000$. The threshold was improved to $n\ge 33$ by Keevash and Mubayi \cite{kemu}.

      \begin{thm}\label{f5}
        For any $k\ge 3$ and sufficiently large $n$, we have $\ex_k^3(n,F_5)=\max\{2\binom{n}{3}, k\ex^3(n,F_5)\}$.
    \end{thm}

    Finally, we consider the case $k$ is fixed but sufficiently large with respect to $\cF$. In the case of graphs, it was shown in \cite{lmz} that if $F$ has a color-critical edge (an edge whose removal decreases the chromatic number) and $k$ is sufficiently large, then for sufficiently large $n$ we have $\ex_k(n,F)=k\ex(n,F)$. On the other hand, if $F$ has no color-critical edge, then this does not hold \cite{chak}.
    We conjecture that a similar characterization holds in the case of hypergraphs.

    \begin{conjecture}
        For any $r$-graph $F$, the following two statements are equivalent.

        (i) There is a $k_0$ such that for every $k\ge k_0$ we have $\ex_k^r(n,F)=k\ex^r(n,F)$ for sufficiently large $n$.

        (ii) There is a constant $c>0$ such that $\ex^r(n,F(1)) \le \ex^r(n,F)-cn^r$ for sufficiently large $n$. 
    \end{conjecture}

    Note that (i) trivially implies (ii). Indeed, if (ii) does not hold, then we take a copy of $K_n^r$ and $k-1$ copies of the extremal graph for $F(1)$. Compared to $k\ex(n,F)$, we gain $\Theta(n^r)$ in the first hypergraph and lose $o(n^r)$ in the other hypergraphs.

    Theorem \ref{f5} shows that the above conjecture holds for $F_5$. We prove another special case. 
       Given a graph $F$, its \textit{$r$-uniform expansion} $F^{r+}$ is obtained by adding $r-2$ new vertices to each edge, such that the $e(F) (r-2)$ new vertices are distinct from each other and are not in $V(F)$.
    An edge of $F$ is said to be \textit{color-critical} if deleting this edge decreases the chromatic number. Given a graph $F$ with chromatic number $t+1$ and a color-critical edge, for $r\le t$ and sufficiently large $n$ we have $\ex^r(n,F^{r+})=|E(T^r(n,t))|$. According to \cite{mv}, this was observed by Alon and Pikhurko, see \cite{gerb} for a proof.

        \begin{thm}\label{thm: sum of edges expansion of color-critical graph}
        Let $r \geq 3$ and let $F$ be a color-critical graph with $\chi(F) = t+1 > r$.
        Then there exists $k_0 = k_0(F, r)$ and $n_0 = n_0(F, r)$ such that for every $k \geq k_0$ and $n \geq n_0$,
        \begin{equation*}
            ex_k^r(n, F^{r+}) = k |E(T^r(n, t))|.
        \end{equation*}
    \end{thm}
    
    In the special case $F$ is a clique, this gives back a result of \cite{zwz}, with an unspecified threshold. 

\section{Proof of Theorems \ref{thm: sum of edges non-r-partite} and \ref{thm: sum of edges non-r-partite large k}}\label{sec2}

We start with a very important lemma that will be used throughout the paper.
For a collection $G_1,\dots,G_k$ of hypergraphs, the \textit{multiplicity} of a hyperedge is the number of hypergraphs among the $G_i$'s that contain it.

\begin{lemma}\label{monot}
    Let $\cF$ be a family of $r$-graphs and $G_1,\dots, G_k$ be a rainbow $\cF$-free collection of $r$-graphs. Then there exists another collection $H_1,\dots,H_k$ satisfying

    (i) every hyperedge has the same multiplicity in both collections.

    (ii) $H_1\supseteq H_2\supseteq \dots \supseteq H_k$.

    (iii) $H_1,\dots,H_k$ is rainbow $\cF$-free.
\end{lemma}

The graph version for a single graph was proved in \cite{Keevash2004} and this version follows from the more general analogous statement in \cite{bks}.
Here we present a different proof.

\begin{proof}
Let us place into $H_i$ each edge with multiplicity at least $i$. Then (i) and (ii) are immediately satisfied, since an edge with multiplicity $j$ is placed into $H_1,H_2,\dots, H_j$. If there is a rainbow copy of a member of $\cF$, its edges are from $H_{i_1},H_{i_2}, \dots H_{i_\ell}$ for some $i_1<i_2<\dots<i_\ell$. In particular, $i_j\ge j$. We go through the hyperedges in this order and each time we pick a new hypergraph from the original collection containing the hyperedge. This is doable, since at step $j$ we have already used $j-1$ hypergraphs, but the hyperedge is from $H_{i_j}$, thus has multiplicity at least $j$. We can pick a new hypergraph always, thus this subgraph was rainbow in the original hypergraph as well, a contradiction completing the proof.
\end{proof}

  \begin{proof}[Proof of Theorem \ref{thm: sum of edges non-r-partite}]
       Let $H_1,\dots, H_k$ be a collection without a rainbow member of $\cF$. By Lemma \ref{monot}, we can assume that $H_1\supseteq H_2\supseteq \dots \supseteq H_k$.
        In particular, $H_\ell,H_{\ell+1},\dots, H_k$ are $F$-free, since a copy of $F$ corresponds to a rainbow copy of $F$. 
        Thus $e(H_{i}) \leq ex^r(n, F)$ for $i \ge \ell$.
        For an $r$-set $e$, denote the multiplicity of $e$ by $m(e) = |\{i \le k: e \in E(H_i)\}|$. 
        Observe that $\sum_{i \le \ell} e(H_i)=\sum_{e \in \binom{n}{r}} \min \{m(e), \ell\}$, since the multiplicity of any $r$-set in the collection $H_1,\dots, H_\ell$ is $\min \{m(e), \ell\}$.
        Then, we have

          \begin{align*}
            \sum_{i \le k} e(H_i)=&\sum_{i \le \ell} e(H_i)+\sum_{\ell+1\le i \le k} e(H_i) \\ \le& \sum_{i \le \ell} e(H_i) +(k-\ell) ex^r(n, {F})\\=& \sum_{e \in \binom{n}{r}} \min \{m(e), \ell\} + (k-\ell) ex^r(n, {F}).
        \end{align*}

        Next, we give an upper bound on $\sum_{e \in \binom{n}{r}} \min \{m(e), \ell\}$.
        \begin{clm}\label{clm: upper bound of sum of multiplicity}
            $\sum_{e \in \binom{n}{r}} \min \{m(e), \ell\} \leq (\ell-1) \binom{n}{r}$.
        \end{clm} 
        \begin{proof}
            Let $F\in\cF$ with $\ell$ hyperedges and $\phi$ be an injective map from $V(F)$ to $V(H_1)$.
            Define the multiplicity of $\phi$ by $m(\phi) = \sum_{e \in {F}} \min \{m(\phi(e)), \ell\}$.
            Then 
            \begin{equation*}
                \sum_\phi m(\phi) = \sum_\phi \sum_{e \in F} \min \{m(\phi(e)), \ell\} = \sum_{e \in F} \sum_\phi \min \{m(\phi(e)), \ell\}.
            \end{equation*}
            Fix a hyperedge $e \in E(F)$.
            For a fixed $r$-set $h$, there are $r! (n-r)_{v-r}$ injective maps such that the image of $e$ is $h$, where $a_b$ denotes $a(a-1)(a-2)\dots (a-b+1)=a!/(a-b)!$.
            Thus $\sum_\phi \min \{m(\phi(e)), \ell\} = r! (n-r)_{v-r} \sum_{h \in \binom{n}{r}} \min \{m(h), \ell\}$.
            Since the size of $F$ is $\ell$, $\sum_{e \in {F}} \sum_\phi \min \{m(\phi(e)), \ell\} = \ell r! (n-r)_{v-r} \sum_{h \in \binom{n}{r}} \min \{m(h), \ell\}$.

            Suppose that the statement does not hold, i.e., $\sum_{e \in \binom{n}{r}} \min \{m(e), \ell\} > (\ell-1) \binom{n}{r}$.
            Then 
            \begin{equation*}
                \sum_\phi m(\phi) = \sum_{e \in F} \sum_\phi \min \{m(\phi(e)), \ell\} > \ell r! (n-r)_{v-r} (\ell-1) \binom{n}{r} = \ell(\ell-1) \frac{n!}{(n-v)!}.
            \end{equation*}
            The total number of injective maps from $V(F)$ to $V(H_1)$ is $n \cdot (n-1) \cdots (n-v+1)$.
 This means that there exists an injective map such that its multiplicity is greater than $\ell (\ell-1)$, i.e., there exists a copy of $F$ such that $\sum_{e \in F} \min \{m(\phi(e)), \ell\} > \ell (\ell-1)$.
            
            Construct an auxiliary bipartite graph with parts $X = E(F)$ and $Y = [k]$.
            Two vertices $e \in X$ and $i \in Y$ are adjacent if and only if $\phi(e)$ is in $H_i$.
            If there is a matching saturating $X$, then we can find a rainbow copy of ${F}$. 
            Otherwise by Hall's condition,
            there exists a  nonempty subset $X'$ of $X$ such that $|\cup_{e \in X'} N(e)| \leq |X'| - 1$, which implies that $m(\phi(e)) \leq |X'| - 1$ for $e \in X'$.
            Hence, we have
            \begin{equation*}
                \sum_{e \in E(F)} \min \{m(\phi(e)), \ell\} \leq |X'| (|X'| - 1) + (\ell - |X'|)\ell \leq \ell(\ell-1),
            \end{equation*}
            a contradiction.
        \end{proof}
      Combining Claim \ref{clm: upper bound of sum of multiplicity} with the inequality above it completes the proof of the Theorem.
    \end{proof}

    \begin{proof}[Proof of Theorem \ref{thm: sum of edges non-r-partite large k}]
        When $\ell = 1$ or $n < v(F)$, it is trivial.
        So we assume that $\ell \geq 2$ and $n \geq v(F)$.
        For any $F$-free $r$-graph $G$, let $\phi:V(F)\to[n]$ be a uniformly random injection. For each edge $e\in E(F)$, we have
        \[
            \Pr(\phi(e)\notin E(G))=\frac{\binom{n}{r}-e(G)}{\binom{n}{r}}.
        \]
        By the union bound, if $e(G)>(1-1/\ell)\binom{n}{r}$, then the probability that at least one edge of $F$ is missing from $G$ is less than $1$. 
        Thus $G$ would contain a copy of $F$, a contradiction. Consequently,
        \[
            ex^r(n,F)\leq\left(1-\frac{1}{\ell}\right)\binom{n}{r},
        \]
        thus $k\geq(\ell-1)\binom{n}{r}-\ell\ex^r(n,F)+\ell\geq\ell$.
        
        The lower bound is given by $k$ copies of an extremal $F$-free hypergraph.
        For the upper bound let $H_1,\dots, H_k$ form a rainbow $F$-free collection with $H_1\supseteq H_2\supseteq \dots, \supseteq H_k$.
        If $|E(H_{\ell})| = ex^r(n, {F})$, then for any $r$-set $h \notin E(H_{\ell})$ there exists a copy of $F$ containing $h$ and other $\ell-1$ hyperedges in $E(H_{\ell})$. 
        If this $r$-set $h$ is a hyperedge, then they form a rainbow copy.
        Since first we can pick a color for $h$ and assign $l-1$ unused colors to remianing $l-1$ hyperedges.
        Thus every $r$-set $h \notin E(H_{\ell})$ cannot be a hyperedge.
        Therefore, $H_1=H_2=\dots=H_\ell\supseteq H_{\ell+1}\supseteq \dots, \supseteq H_k$ and we are done.
   
        If $|E(H_{l})| \leq ex^r(n, {F}) - 1$, then we have

            \begin{align*}
            \sum_{i \le k} |E(H_i)|=&\sum_{i \le \ell} |E(H_i)|+\sum_{\ell+1\le i \le k} |E(H_i)| \\ 
            \le& \sum_{i \le \ell} |E(H_i)| +(k-\ell) (ex^r(n, {F}) -1)\\
            \leq& (l-1) \binom{n}{r} + (k-\ell) (ex^r(n, {F}) - 1)\\
            \leq& k ex^r(n, {F}).
        \end{align*}
        Here the first inequality is trivial, the second inequality uses Claim \ref{clm: upper bound of sum of multiplicity}, and the last inequality uses the assumptions of the theorem.
    \end{proof}

\section{Proof of Theorem \ref{thm: sum of edges r-partite}}\label{sec3}

    Let $H_1,\dots, H_k$ be a rainbow $F$-free collection with $H_1\supseteq H_2\supseteq \dots, \supseteq H_k$.
    Let $v = |V(F)|$.
   
    It is easy to see that $H_\ell$ is $F$-free.
    Otherwise, each hyperedge of that copy of $F$ is contained in at least $\ell$ hypergraphs $H_i$, thus we can find a rainbow copy of $\mathcal{F}$ greedily.
    Let
    \begin{equation*}
        H_{\ell-1}' = K_n^r \setminus H_{\ell-1}.
    \end{equation*}
    
    Since $H_\ell$ is $F$-free, by a result of Erd\H{o}s \cite{Erdos1964}, we have
    \begin{equation*}
        |E(H_\ell)| = o(n^r).
    \end{equation*}
    Suppose that $\sum_{i=1}^{k} |E(H_i)| > (\ell-1) \binom{n}{r}$. This implies
    \begin{equation*}
        |E(H'_{\ell-1})| = o(n^r).
    \end{equation*}

    In the proofs, we consider copies of $F$ on the vertex set of our collection, independently of which hypergraphs $H_i$ contain their hyperedges. In particular, we consider those $r$-sets that are not hyperedges of any $H_i$.
    
    \begin{lemma}\label{lem: sets with large degree in Hl-1'}
        There is a constant $c$ such that for sufficiently large $n$ and each hyperedge $e$ of $H_\ell$, there exists a proper subset $S \subseteq e$ such that
        \begin{equation*}
            d_{H_{\ell-1}'}(S) \geq c n^{r-|S|}.
        \end{equation*} 
    \end{lemma}
    \begin{proof}    
        For each hyperedge $e$ in $H_\ell$, there are $(n-r)_{v-r}  \geq \frac{1}{2} n^{v-r}$ copies of $F$ containing $e$ for sufficiently large $n$.
        Each such copy contains at least one hyperedge in $H_{\ell-1}'$.
        Suppose otherwise, we can choose $e$ from $H_\ell$, the remaining $\ell-1$ hyperedges appear in $H_{\ell-1}$ and thus in $H_i$ where $i \le \ell-1$.
        We can choose these $\ell-1$ hyperedges from different hypergraphs and thus get a rainbow copy of $F$, a contradiction.
        
        Consider a hyperedge $h$ in $H_{\ell-1}'$.
        If $h$ does not intersect with $e$, then the number of copies of $F$ that contain both $h$ and $e$ is at most $c_1 n^{v-2r}$.
        Thus, there are at most $|E(H_{\ell-1}')| c_1 n^{v-2r}$ copies of $F$ that contain $e$ and a hyperedge $h$ of $H_{\ell-1}'$ where $h \cap e = \emptyset$.
        Since $|E(H_{\ell-1}')| = o(n^r)$, we have $|E(H_{\ell-1}')| c_1 n^{v-2r} \leq \frac{1}{4} n^{v-r}$ for sufficiently large $n$.
        
        Therefore, there are at least $\frac{1}{4} n^{v-r}$ copies of $F$ containing $e$ and at least one hyperedge in $H_{\ell-1}'$ that intersects with $e$.
        For such a hyperedge $h$ with $h \cap e = T \neq \emptyset$, there are at most $c_2 n^{v-2r+|T|}$ copies of $F$ that contain both $e$ and $h$.
        Let $A_i(e) = \{h \in H_{\ell-1}': |e \cap h| = i\}$, where $i \le r-1$.
        Then we have $\sum_{i=1}^{r-1} |A_i(e)| c_2 n^{v-2r+i} \geq \frac{1}{4} n^{v-r}$.
        By pigeon-hole principle, for each hyperedge $e$ in $H_\ell$, there exists an $i$ such that $|A_i(e)| \geq \frac{n^{r-i}}{4 c_2 (r-1)}$.
        Each hyperedge contains $\binom{r}{i}$ $i$-sets.
        By pigeon-hole principle, there exists a proper subset $S \subseteq e$ with size $i$ such that the number of hyperedges in $H_{\ell-1}'$ containing $S$ is at least $\frac{n^{r-i}}{4 c_2 (r-1) \binom{r}{i}}$.
        Let $c = \min_{1\leq i\leq r-1}\frac{1}{4 c_2 (r-1) \binom{r}{i}}$.
        Then this proof is completed.
    \end{proof}

    
    Let $0 < c_1 < c_2 < \dots < c_{r-1} < c$ be constants. We choose them in decreasing order of their indices so that $c_i$ is small enough with respect to $c_{i+1}$ (a precise dependence will be given later).
    We say that a set $S$ in $H_{\ell-1}'$ is \textit{good} if it is a minimal nonempty set with $d_{H_{\ell-1}'}(S) \geq c_{|S|} n^{r-|S|}$, i.e., for any nonempty proper subset $T \subseteq S$ has $d_{H_{\ell-1}'}(T) < c_{|T|} n^{r-|T|}$ where $c_{|T|} < c_{|S|}$.
    Let $U_i$ be the set of good $i$-sets, $i \in [r-1]$. 

    \begin{lemma}\label{lem: sets with small degree in Hl}
        For $1 \leq i \leq r-1$ and any $\epsilon > 0$, let $S$ be an $i$-set in $U_i$.
        Then $d_{H_\ell}(S) \leq \epsilon n^{r-i}$ for sufficiently large $n$.
    \end{lemma}
    \begin{proof}
        Suppose otherwise, there exists a set $S$ in $U_i$ with $d_{H_\ell}(S) \geq \epsilon n^{r-i}$ for some $\epsilon$ and sufficiently large $n$.
        Let $L(S, H_\ell) = \{T: T \cup S \in E(H_\ell)\}$ be the link hypergraph of $S$.
        Thus the size of $L(S, H_\ell)$ is at least $\epsilon n^{r-i}$. Recall that $F$ is an $r$-partite $r$-uniform hypergraph.
        Choose a non-isolated set $S' \subseteq V(F)$ with size $i$, i.e., a set of size $i$ with vertices from distinct parts of $F$ that is contained in at least one hyperedge.
        Then $L(S', F)$ is an $(r - i)$-partite $(r - i)$-uniform hypergraph.

        We will embed the hyperedges containing $S'$ in $F$ to hyperedges containing $S$ in $H_\ell$.
        For $i \leq r-2$, we apply Erd\H{o}s's theorem \cite{Erdos1964} to show that the extremal number of $L(S', F)$ is $o(n^{r-i})$. Since there are $\Omega(n^{r-i})$ hyperedges in $L(S, H_\ell)$,
        we can embed $S'$ to $S$ and embed hyperedges containing $S'$ to hyperedges containing $S$.
        For $i = r - 1$, the link hypergraph $L(S', F)$ is a set of isolated vertices and the size of $L(S, H_\ell)$ is $\Omega(n)$.
        We can also embed $S'$ to $S$ and embed hyperedges containing $S'$ to hyperedges containing $S$.
        
        Moreover, by the well-known supersaturation phenomenon \cite{ersim}, we find not only one, but $c_F' n^{v(L(S', F))}$ copies of $L(S', F)$ in $L(S, H_l)$. 
        We extend these copies of $L(S', F)$ to copies of $F$ (not necessarily in $H_1$) by embedding the remaining vertices of $F$.
        Then there are $c_F n^{v - i}$ copies of $F$ containing $S$.
        Each such copy contains a hyperedge $h$ in $H_{\ell-1}'$.
        By construction, every edge of $F$ containing $S'$ was embedded into $H_{\ell}$.
        Thus an edge $h$ of this copy that lies in $H_{\ell-1}'$ cannot contain all of $S$, i.e., $h \cap S$ is a proper subset of $S$.
        If $h \cap S = T \neq \emptyset$, then by the definition of a good set, $d_{H_{\ell-1}'}(T) < c_{|T|} n^{r-|T|}$.
        Obviously, the number of copies of $F$ containing $h$ and $S$ is at most $O(n^{v - i - r + |T|})$.
        If $h \cap S = \emptyset$, then the number of copies of $F$ containing $h$ and $S$ is at most $O(n^{v - i - r})$.
        Thus, the total number of such copies of $F$ that contain $S$ is at most
        \begin{equation*}
            \sum_{j = 1}^{i-1} \binom{i}{j} c_{j} n^{r-j} O(n^{v - i - r + j}) + |E(H_{\ell-1}')| O(n^{v - i - r}).
        \end{equation*}
        We choose each $c_{j}$ small enough in decreasing order of their indices such that this value is less than $c_F n^{v - i}$ for sufficiently large $n$. In this way, we obtain
        a contradiction completing the proof.
    \end{proof}

    Choose $\epsilon > 0$ such that $\epsilon (k-\ell+1) (2^r-2) < 1$ and $0< c_1 < c_2 \dots < c_{r-1} < c$ small enough in decreasing order of their indices, so that Lemma \ref{lem: sets with small degree in Hl} can be applied with parameter $\epsilon c_i$ for every $i$. 
    Since
    \begin{equation*}
        |U_i| c_i n^{r-i} \leq \sum_{S \in U_i} d_{H_{\ell-1}'} (S) \leq \sum_{S \in \binom{V(H)}{i}} d_{H_{\ell-1}'} (S) = \binom{r}{i} |E(H_{\ell-1}')|,
    \end{equation*}
    we have
    \begin{equation*}
        |U_i| \leq \binom{r}{i} \frac{|E(H_{\ell-1}')|}{c_i n^{r-i}}.
    \end{equation*}
    By Lemma \ref{lem: sets with large degree in Hl-1'}, every hyperedge $e$ in $H_\ell$ contains a good set.
    Thus
    \begin{equation*}
        |E(H_\ell)| \leq \sum_{i=1}^{r-1} \sum_{S \in U_i} d_{H_\ell}(S) \leq \sum_{i=1}^{r-1} |U_i| \epsilon c_i n^{r-i} \leq \epsilon \sum_{i=1}^{r-1} \binom{r}{i} \frac{|E(H_{\ell-1}')|}{c_i n^{r-i}} c_i n^{r-i} = \epsilon |E(H_{\ell-1}')| (2^r - 2).
    \end{equation*}
    Hence we have
    \begin{align*}
        \sum_{i=1}^{k} |E(H_i)| =& \sum_{i=1}^{\ell-1} |E(H_i)| + \sum_{i=\ell}^{k} |E(H_i)| \\
                                        \leq& (\ell-1) \binom{n}{r} - |E(H_{\ell-1}')| + (k-\ell+1) \epsilon |E(H_{\ell-1}')| (2^r - 2) \\
                                        =& (\ell-1) \binom{n}{r} - (1- \epsilon (k-\ell+1) (2^r-2)) |E(H_{\ell-1}')| \\
                                        <& (\ell-1) \binom{n}{r},
    \end{align*}
    where the last inequality follows from the choice of $\epsilon$.
   This is a contradiction, thus the proof of Theorem \ref{thm: sum of edges r-partite} is complete.

\section{Proof of Theorem \ref{f5}}


We will adapt the proof that determines the Tur\'an number of $F_5$ for sufficiently large $n$ from \cite{ff}.

\begin{proof}[Proof of Theorem \ref{f5}] Let $H_1,\dots, H_k$ be a rainbow $F_5$-free collection with $H_1\supseteq H_2\supseteq \dots \supseteq H_k$.
We are going to use the progressive induction method of Simonovits \cite{sim}. In our particular case it works the following way. Let us denote the claimed bound by $f(n)$. We are going to show that for $n> n_0$, either $H_1$ is $F_5$-free, thus each $H_i$ is $F_5$-free and we are done, or there is a vertex incident to less than $f(n)-f(n-1)$ hyperedges (in the whole collection). The total number of hyperedges in a rainbow $F$-free collection on $n_0$ vertices is at most $f(n_0)+c$ for some constant $c$, since $n_0$ is a constant. Then on $n_0+1$ vertices we have at most $f(n_0+1)+c-1$ hyperedges since the deleted vertex was in at most $f(n_0+1)-f(n_0)-1$ hyperedges (or $H_1$ is $F_5$-free and we are done). Continuing this way, we obtain that on $n_0+c$ vertices we have at most $f(n_0+c)$ hyperedges. Observe that $f(n)-f(n-1)=\max\{n^2, kn^2/9\}+O(n)$. We will show a vertex that is incident to at most $\alpha n^2$ hyperedges for some $\alpha<\max\{1,k/9\}$. Since $n$ is sufficiently large, this is smaller than $f(n)-f(n-1)$.

    Let $G$ denote the shadow graph of $H_1$, i.e., the graph consisting of edges that are contained in at least one hyperedge of $H_1$. If $G$ is $K_4$-free then $H_1$ is $F_5$-free, thus each $H_i$ is $F_5$-free and we are done. 

Let us fix a constant $t \geq k + 2$ and let $s=\min \{t,\omega(G)\}$, where $\omega(G)$ is the order of the largest clique of $G$. Let us fix a $K_s$ with vertex set $S$. For each of the $\binom{s}{2}$ edges of this clique, we pick a vertex that extends it to a hyperedge of $H_1$, let $T$ denote the set of these at most $\binom{s}{2}$ vertices. Note that $T$ may intersect $S$. Let $X=V(G)\setminus (S\cup T)$. Let $G'$ and $G''$ denote the shadow graphs of $H_2$ and $H_3$ restricted to $X$.
Since $|S \cup T| = O(1)$, the total number of hyperedges intersecting $S \cup T$ in at least two vertices is $O(n)$.

\begin{clm}\label{clam}
    For each edge $uv$ of $G''$, there is at most one vertex of $S$ that extends $uv$ to a hyperedge of $H_3$. Moreover, if there is such a vertex $w$, then $w$ is the only vertex of $S$ that extends $uv$ to a hyperedge of $H_2$.
\end{clm}

\begin{proof}[Proof of Claim]
Assume that $uv$ is extended by both $x$ and $y$ from $S$. By definition, there is a hyperedge in $H_1$ containing $x,y$ and a vertex $z$ not in $X$. Thus $u,v,x,y,z$ are 5 vertices with $\{u,v,x\}$ in $H_3$, $\{u,v,y\}$ in $H_3$ (thus also in $H_2$) and $\{x,y,z\}$ in $H_1$. They form a rainbow $F_5$, a contradiction. Assume now $w \neq w'$ and $\{u,v,w\}$ is in $H_3$ and $\{u,v,w'\}$ is in $H_2$. There is a hyperedge in $H_1$ containing $w,w'$ and a vertex $z$ not in $X$. They together form a rainbow $F_5$, a contradiction.
\end{proof}

Let us return to the proof of the Theorem. 

\textbf{Case 1.} We have $s=\omega(G)$. Then for each edge $uv$ of $G$, there are at most $s-2$ vertices in $S$ that extend $uv$ to a hyperedge of $H_1$, since they from a clique in $G$ with $u,v$. By Claim \ref{clam}, at most $\max\{s-2,k-1\}$ vertices extend $uv$ to a hyperedge of some $H_i$ with $i>1$. In the following, we also say that a vertex extends $uv$ in $H_i$ instead of saying that it extends $uv$ to a hyperedge of $H_i$.
If there exists a vertex of $S$ that extends $uv$ in $H_3$, then by Claim \ref{clam} it is also the only vertex of $S$ that extends $uv$ in $H_2$.
Since $H_1\supseteq H_2\supseteq \dots \supseteq H_k$, it is the only vertex that can extend $uv$ in $H_i$ for $i > 2$.
If there is no vertex of $S$ that extends $uv$ in $H_3$, then there are at most $s-2$ vertices of $S$ that extend $uv$ in $H_2$.

Therefore, the total number of hyperedges containing exactly one vertex from $S$ and two from $X$ is at most $\max\{2s-4,k+s-3\}|E(G)|\le \max\{2s-4,k+s-3\}(s-1)n^2/2s$ using Tur\'an's theorem \cite{Tu}. This implies that there is a vertex $v$ of $S$ that is incident to at most $\max\{2s-4,k+s-3\}(s-1)n^2/2s^2 +O(n)$ hyperedges in total. If the maximum is given by $2s-4$, then this is at most $(2s-4)(s-1)n^2/2s^2 + O(n)<(2s^2-5s)n^2/2s^2+O(n)< f(n)-f(n-1)$. Here we used the fact that $s$ is bounded above by a constant $t$. 

If the maximum is given by $k+s-3$, then $v$ is incident to at most $(k+s-3)(s-1)n^2/2s^2 + O(n)$ hyperedges. One can check that $(k+s-3)(s-1)/2s^2<1$ for $k\le 9$, thus this upper bound is less than $f(n)-f(n-1)$ for $k\le 9$. If we increase $k$ by 1, then this increases by $(s-1)/2s^2$, while $f(n)-f(n-1)$ increases by $n^2/9+O(n)$. Since $s\ge 4$, $f(n)-f(n-1)$ increases more, thus the degree of $v$ is less than $f(n)-f(n-1)$ for all $k$.

\textbf{Case 2.} We have $s=t$, i.e., there are cliques of order more than $s$ in $G$. Then we lose two statements from the argument in Case 1: for an edge $uv$ we can have more than $s-2$ vertices in $S$ that extend it to a hyperedge of $H_1$ (and similarly for $H_2$), and we do not have the same upper bound on $|E(G)|$. 



Recall that we have not considered the hyperedges containing more than one vertex in $S\cup T$. Let us now consider those as well. We say that an edge of $G$ is \textit{heavy} if it is contained in at least three hyperedges of $H_1$. We show that if there are only heavy edges inside $S$, then Claim \ref{clam} remains true.

\begin{clm}\label{clam2} Assume that each edge inside $S$ is heavy.
    Then for each edge $uv$ of $G$, there is at most one vertex of $S$ that extends them to a hyperedge of $H_3$. Moreover, if there is such a vertex $w$, then $w$ is the only vertex of $S$ that extends $uv$ to a hyperedge of $H_2$.
\end{clm}

\begin{proof}[Proof of Claim]
Assume that $uv$ is extended by both $x$ and $y$ from $S$ to a hyperedge of $H_3$. By definition, there is a hyperedge in $H_1$ containing $x,y$ and a vertex $z$ distinct from $u$ and $v$. Then these three hyperedges form a rainbow $F_5$, a contradiction. Assume now that $\{u,v,w\}$ is in $H_3$ and $\{u,v,w'\}$ is in $H_2$. There is a hyperedge in $H_1$ containing $w,w'$ and a vertex $z$ distinct from $u$ and $v$. Together, these three hyperedges form a rainbow $F_5$, a contradiction.
\end{proof}

By the above claim, if each edge inside $S$ is heavy, then for each edge $uv$ that contains $i$ vertices from $S$, there are at most $\max \{2s-2i,s+k-i-1\}=2s-2i$ hyperedges containing $uv$ and a vertex from $S$.
Observe that this is sharp for the construction of two complete 3-graphs and $k-2$ empty 3-graphs. 
So the total degree sum over $S$ is at most $\binom{s}{0} \binom{n-s}{2} 2s + \binom{s}{1} \binom{n-s}{1} (2s-2) + \binom{s}{2} \binom{n-s}{0} (2s-4) = s(n-1)(n-2)$.
This implies that the average degree in $S$ is at most $f(n)-f(n-1)$. We are done unless we have equality for each vertex of $S$. In this case every 3-set intersecting $S$ is a hyperedge of both $H_1$ and $H_2$. Then, it is easy to see that any hyperedge in $H_3$ would create a rainbow $F_5$. Therefore, $H_3,\dots, H_k$ are empty and we are done. 

It is left to deal with the case that there is a non-heavy edge inside $S$, and analogously inside every $K_s$. This means that the graph of heavy edges is $K_s$-free. 
The number of heavy edges in $X$ is at most $ex^2(|X|, K_s) = (1-\frac{1}{s-1} + o(1))\frac{|X|^2}{2}$ by Tur\'an's theorem \cite{Tu}.
If a pair of vertices in $X$ is not an edge in $G$, then there is no hyperedge that contains this pair.
If a pair of vetices in $X$ is non-heavy edge, then there are at most $2k$ hyperedges containing this pair.
Observe that Claim \ref{clam} is still true for $u,v\in X$. Therefore, the total number of hyperedges containing a vertex from $S$ and $u,v$ from $X$ is at most $\max\{2s,k+s-1\}=2s$. 
Since $s = t > k$, $2s x + 2k (\binom{|X|}{2} - x)$ is increasing in $x$.
This implies that the total degree sum over $S$ is at most $2s ex^2(|X|, K_s) + 2k (\binom{|X|}{2} - ex^2(|X|, K_s)) + O(n) = (s-1 + \frac{k-1}{s-1} + o(1))n^2 + O(n)$.
Therefore, there is a vertex of $S$ with total degree less than $2\binom{n-1}{2}$, completing the proof.
\end{proof}

\section{Proof of Theorem \ref{thm: sum of edges expansion of color-critical graph}}





   We will use the following result of Tang, Xu and Yan \cite{TangXuYan}.

    \begin{thm}[\cite{TangXuYan}, Prop. 2.3]\label{thm: stability of expansion}
        Fix $t \ge r \geq 2$ and a graph $F$ with $\chi(F) = t+1$.
        For every $\epsilon_1 > 0$, there exist $n_1 = n_1(r, t, \epsilon_1) > 0$, $\eta = \eta(r, t, \epsilon_1) > 0$, such that if $n > n_1$ and $H$ is an $n$-vertex $F^{r+}$-free $r$-uniform hypergraph with $|E(H)| \geq |E(T^r(n, t))| - \eta \binom{n}{r}$, then one can add or remove at most $\epsilon_1 n^r$ hyperedges from $H$ to make it isomorphic to $T^r(n, t)$.
    \end{thm}

    Note that in the case $F=K_{t+1}$, this was proved by Mubayi \cite{mub}. It implies the following bound on the Tur\'an number of $F^{r+}$. Note that the following result was claimed in \cite{mv} without a detailed proof, see \cite{pttw} for a proof.
    
    \begin{thm}[\cite{mv}]\label{thm: extremal number of expansion}
        Fix integers $t \geq r \geq 2$ and a graph $F$ with $\chi(F) = t+1$.
        Then
        \begin{equation*}
            ex_r (n, F^{r+}) =(1+o(1)) |E(T^r (n, t))|.
        \end{equation*}
        In particular, for every $\epsilon_2 > 0$ and sufficiently large $n$, we have $ex_r (n, F^{r+}) \leq |E(T^r (n, t))| + \epsilon_2 n^r$.
    \end{thm}

Now we are ready to prove Theorem \ref{thm: sum of edges expansion of color-critical graph}. Recall that the case of cliques was proved in \cite{zwz}. Our proof is an extension of their proof.

    \begin{proof}[Proof of Theorem \ref{thm: sum of edges expansion of color-critical graph}]The lower bound is given by $k$ copies of $T^r(n, t)$. For the upper bound,
        let $H_1,\dots, H_k$ be a collection of hypergraphs without a rainbow copy of $F^{r+}$ that maximizes the total number of hyperedges. We denote this collection by $H$ and denote the total number of hyperedges in a collection $H$ by $e(H)$.

        Choose $n_0$ large enough so that $\frac{n_0}{2}$ satisfies the conditions of Theorem \ref{thm: stability of expansion} and Theorem \ref{thm: extremal number of expansion}.
        We prove the upper bound by conducting a contradiction.
        Assume that $e(H)>k |E(T^r(n, t))|$, and let $f(H)=e(H)-k |E(T^r(n, t))|$.
        By \ref{monot}, we assume that $H_k \subseteq H_{k-1} \subseteq \cdots \subseteq H_1$.
        For convenience, let $v_F = v(F^{r+})$, 
        $\ell = e(F^{r+})$ and $c_0 = \frac{\ell - 1}{k}$.
        \begin{clm}\label{clm: subgraph with large minimum degree}
            If $k \geq k_0$, $N \geq n_0$ and $f(H) > 0$, then there exists a collection $H'_1,\dots, H'_k$ (denoted by $H'$) on $N' > \frac{N}{2}$ vertices such that each $H'_i$ is a subhypergraph of $H_i$, $f(H') > 0$ and
            \begin{equation*}
                d_{H'}\left(v\right) \geq k\left(|E\left(T^r\left(N', t\right)\right)| - |E\left(T^r\left(N'-1, t\right)\right)| - c_1 \left(N'\right)^{r-1}\right)
            \end{equation*}
            for all $v \in V\left({H}'\right)$, where $c_1 \geq 3 \epsilon_2$.
        \end{clm}
        \begin{proof}[Proof of Claim \ref{clm: subgraph with large minimum degree}]
            Denote $H$ by $G_0$.
            If there exists a vertex $v_1$ with
            \begin{equation*}
                d_{G_0}\left(v_1\right) < k\left(|E\left(T^r\left(N, t\right)\right)| - |E\left(T^r\left(N-1, t\right)\right)| - c_1 N^{r-1}\right),
            \end{equation*}
            then delete $v_1$ and denote the new collection of hypergraphs by $G_1$.
            This operation will not create a rainbow copy of $F^{r+}$ and
            \begin{align*}
                f\left(G_1\right) - f\left(G_0\right) =& \left(e\left(G_1\right) - k |E\left(T^r\left(N-1, t\right)\right)|\right) - \left(e\left(G_0\right) - k |E\left(T^r\left(N, t\right)\right)|\right) \\
                =& k \left(|E\left(T^r\left(N, t\right)\right)| - |E\left(T^r\left(N-1, t\right)\right)|\right) - d_{G_0}\left(v_1\right) \\
                >& k c_1 N^{r-1}.
            \end{align*}
            Suppose that this operation continues until there exists a vertex $v_{\left\lceil\frac{N}{2}\right\rceil}$ in $G_{\left\lceil\frac{N}{2}\right\rceil-1}$ with $d_{G_{\left\lceil\frac{N}{2}\right\rceil-1}}\left(v_{\left\lceil\frac{N}{2}\right\rceil}\right) < k\left(|E\left(T^r\left(\left\lfloor\frac{N}{2}\right\rfloor+1, t\right)\right)| - |E\left(T^r\left(\left\lfloor\frac{N}{2}\right\rfloor, t\right)\right)| - c_1 \left(\left\lfloor\frac{N}{2}\right\rfloor+1\right)^{r-1}\right)$.
            Delete $v_{\left\lceil\frac{N}{2}\right\rceil}$ and denote the new collection of hypergraphs by $G_{\left\lceil\frac{N}{2}\right\rceil}$.
            So we have
            \begin{align*}
                f\left(G_{\left\lceil\frac{N}{2}\right\rceil}\right) - f\left(G_{\left\lceil\frac{N}{2}\right\rceil-1}\right) =& \left(e\left(G_{\left\lceil\frac{N}{2}\right\rceil}\right) - k |E\left(T^r\left(\left\lfloor\frac{N}{2}\right\rfloor, t\right)\right)|\right) \\
                &- \left(e\left(G_{\left\lceil\frac{N}{2}\right\rceil-1}\right) - k \left|E\left(T^r\left(\left\lfloor\frac{N}{2}\right\rfloor+1, t\right)\right)\right|\right) \\
                =& k\left(|E\left(T^r\left(\left\lfloor\frac{N}{2}\right\rfloor+1, t\right)\right)| - |E\left(T^r\left(\left\lfloor\frac{N}{2}\right\rfloor, t\right)\right)|\right) - d_{G_{\left\lceil\frac{N}{2}\right\rceil-1}}\left(v_{\left\lceil\frac{N}{2}\right\rceil}\right) \\
                >& k c_1 \left(\left\lfloor\frac{N}{2}\right\rfloor+1\right)^{r-1}.
            \end{align*}
            Thus we have
            \begin{align*}
                f\left(G_{\left\lceil\frac{N}{2}\right\rceil}\right) >& k c_1 \left(\left(\left\lfloor\frac{N}{2}\right\rfloor+1\right)^{r-1} + \left(\left\lfloor\frac{N}{2}\right\rfloor+2\right)^{r-1} + \cdots + N^{r-1}\right) + f\left(G_0\right) \\
                \geq& k c_1 \left\lceil\frac{N}{2}\right\rceil \left(\left\lfloor\frac{N}{2}\right\rfloor+1\right)^{r-1} + f\left(G_0\right) \\
                >& k c_1 \left\lceil\frac{N}{2}\right\rceil \left(\left\lfloor\frac{N}{2}\right\rfloor+1\right)^{r-1}.
            \end{align*}
            On the other hand, let $G_{\left\lceil\frac{N}{2}\right\rceil, i}$ be the $i$-th member of $G_{\left\lceil\frac{N}{2}\right\rceil}$ for $i\le k$. 
            Then
            \begin{align*}
                f\left(G_{\left\lceil\frac{N}{2}\right\rceil}\right) \leq& k \left|E\left(G_{\left\lceil\frac{N}{2}\right\rceil, \ell}\right)\right| + \left(\ell - 1\right) \left(\left|E\left(G_{\left\lceil\frac{N}{2}\right\rceil, 1}\right)\right| - |E\left(G_{\left\lceil\frac{N}{2}\right\rceil, \ell}\right)\right)| - k \left|E\left(T^r\left(\left\lfloor\frac{N}{2}\right\rfloor, t\right)\right)\right| \\
                \leq& k \left(\left|E\left(T^r\left(\left\lfloor\frac{N}{2}\right\rfloor, t\right)\right)\right| + \epsilon_2 \left\lfloor\frac{N}{2}\right\rfloor^r\right) + \left(\ell - 1\right) \binom{\left\lfloor\frac{N}{2}\right\rfloor}{r} - k \left|E\left(T^r\left(\left\lfloor\frac{N}{2}\right\rfloor, t\right)\right)\right| \\
                <& k \epsilon_2 \left\lfloor\frac{N}{2}\right\rfloor^r + \left(\ell - 1\right) \left\lfloor\frac{N}{2}\right\rfloor^r \leq 2 k \epsilon_2 \left\lfloor\frac{N}{2}\right\rfloor^r.
            \end{align*}
            The second inequality follows from Theorem \ref{thm: extremal number of expansion}, and the last inequality follows since we can choose $k$ large enough such that $\epsilon_2 \geq c_0$.
            Since $c_1 \geq 3 \epsilon_2$, we have
            \begin{equation*}
                k c_1 \left\lceil\frac{N}{2}\right\rceil \left(\left\lfloor\frac{N}{2}\right\rfloor+1\right)^{r-1} > 2 k \epsilon_2 \left\lfloor\frac{N}{2}\right\rfloor^r,
            \end{equation*}
            a contradiction.
            Thus, this operation stops with $n' > \frac{N}{2}$, and the resulting hypergraph satisfies $f\left(H'\right) > 0$ and $d_{H'}\left(v\right) \geq k\left(|E\left(T^r\left(N', t\right)\right)| - |E\left(T^r\left(N'-1, t\right)\right)| - c_1 \left(N'\right)^{r-1}\right)$ for every $v \in V\left(H'\right)$.
        \end{proof}
        By Claim \ref{clm: subgraph with large minimum degree}, we can obtain a subcollection $H'$ with large minimum degree.
        For convenience, denote the number of vertices of $H'$ by $n$. 
        It is easy to see that $H'_k \subseteq H'_{k-1} \subseteq \cdots \subseteq H'_1$.
        Since $f(H') > 0$ and $e(H') \leq k |E(H'_{\ell})| + (\ell-1) \binom{n}{r}$, we have
        \begin{equation*}
            |E(H'_{\ell})| > |E(T^r(n, t))| - c_0 \binom{n}{r}.
        \end{equation*}
        Since $k (|E(T^r(n, t))| - |E(T^r(n-1, t))| - c_1 n^{r-1}) \leq d_{H'} (v) \leq k d_{H'_{\ell}} (v) + (\ell - 1) \binom{n-1}{r-1}$, we have
        \begin{align*}
            d_{H'_{\ell}} (v) \geq& |E(T^r(n, t))| - |E(T^r(n-1, t))| - c_1 n^{r-1} - c_0 \binom{n-1}{r-1} \\
            \geq& |E(T^r(n, t))| - |E(T^r(n-1, t))| - c_1' n^{r-1},
        \end{align*}
        where $c_1' = c_1 + \frac{c_0}{(r-1)!}$.
        
        Applying Theorem \ref{thm: stability of expansion} to $H_{\ell}'$ with $c_0 \leq \eta$, there exists a $t$-partition with at least $|E(T^r(n, t))| - \epsilon_1 n^r$ hyperedges which intersect each part in at most one vertex.
        Among all such partitions $P = (V_1, V_2, \dots, V_{t})$, choose one that maximizes $|H_{\ell}' \cap K(P)^r|$, where $K(P)^r$ is the complete $t$-partite $r$-uniform hypergraph with partition $P$.
        Let $H_{\ell}'' = K(P)^r \setminus H_{\ell}'$.
        By Theorem \ref{thm: stability of expansion}, we have
        \begin{equation*}
            |E(H_{\ell}'')| \leq \epsilon_1 n^r \text{ and } \left||V_i| - \frac{n}{t}\right| \leq c_2 n \text{ for all } i \le t,
        \end{equation*} 
        where $c_2 \leq \frac{1}{4t}$ and satisfies $\binom{t-2}{r-1} \left(\frac{1}{t} + c_2\right)^{r-1} \leq \binom{t-2}{r-1} \left(\frac{1}{t}\right)^{r-1} + \frac{1}{16} \binom{t-2}{r-2}\left(\frac{1}{t}\right)^{r-1}$. In particular, each part $V_i$ has order at least $3n/4t$.
        
        Let $X = \{v \in V(H'): d_{H_{\ell}''}(v) \geq c_3 n^{r-1}\}$, 
        where $c_3 > \left(\frac{r^2 \epsilon_1}{2(1-c_0)(r-2)!}\right)^{\frac{1}{2}}$.
        Since $\sum_{v \in V(H_{\ell}'')} d_{H_{\ell}''}(v) = r e(H_{\ell}'')$, we have $|X| \leq r \frac{\epsilon_1 n^r}{c_3 n^{r-1}} = \frac{r \epsilon_1 n}{c_3}$.
        
        We call a pair of vertices $(u, v)$ \textit{heavy} if there are at least $(r-2) \ell + v_F + 1$ hyperedges in $H_{\ell}'$ such that they contain $(u, v)$ and any two of them share only $(u, v)$.
        For $u \in V(H')$ and $i \le t$, let $N_i(u) = \{v \in V_i \setminus (X \cup \{u\}): (u, v) \text{ is heavy}\}$.
        We say that $V_i$ is \textit{rich} for $u$ if $|N_i(u)| \geq c_4 n$ where $\max\left\{\frac{r \epsilon_1}{c_3}, v_F c_5\right\} < c_4 < \min\left\{\frac{3}{4t} - \frac{r\epsilon_1}{c_3} - c_5, \left(\frac{15}{16}\binom{t-2}{r-2}t^{-(r-1)} - c_1' - \frac{r\epsilon_1}{c_3}\right)/2\right\}$, where $c_5 = 2(r-1)c_3 \binom{t-2}{r-2}^{-1} \left(\frac{3}{4t}\right)^{-r+2}$. Otherwise we say that $V_i$ is \textit{poor} for $u$.
        
        \begin{clm}\label{clm: vertex outside X has many heavy partner}
            Let $u \in V_i \setminus X$.
            Then for every $j \neq i$, all but at most $c_5 n$  vertices in $V_j$ form a heavy pair with $u$ and $V_j$ is rich for $u$.
        \end{clm}
        \begin{proof}[Proof of Claim \ref{clm: vertex outside X has many heavy partner}]
            For each pair $(u, v)$ with $v \in V_j, i \neq j$, there are at least
            \begin{equation*}
                \binom{t-2}{r-2} \left(\frac{3n}{4t}\right)^{r-2} = c_6 n^{r-2}
            \end{equation*}
            $r$-sets that contain $u,v$ and intersect each part in at most one vertex, i.e., hyperedges in $K(P)^r$.
            If $(u, v)$ is not a heavy pair, then there are at most 
            \begin{equation*}
                (r-2) ((r-2) \ell + v_F) n^{r-3}
            \end{equation*}
            hyperedges in $H_{\ell}'$ containing them.
            Indeed, consider the link hypergraph $L((u, v), H_{\ell}')$ of $(u, v)$ in $H_{\ell}'$.
            The matching number of $L((u, v), H_{\ell}')$ is less than $(r-2) \ell + v_F + 1$.
            Choose a maximum matching, then all other hyperedges intersect this matching.
            Thus the total number of hyperedges in this link hypergraph is at most $(r-2) ((r-2) \ell + v_F) n^{r-3}$. This implies that the number of hyperedges in $H_\ell''$ containing $u,v$ is at least $c_6 n^{r-2}-(r-2) ((r-2) \ell + v_F) n^{r-3}\ge c_6 n^{r-2}/2$.
            
            For each hyperedge containing $u$ in $H_{\ell}''$, there are $r-1$ possible choices for a vertex $v$ that does not form a heavy pair with $u$.
            Thus the number of vertices outside $V_i$ which do not form a heavy pair with $u$ is at most 
            \begin{equation*}
                \frac{2(r-1) d_{H_{\ell}''}(u)}{c_6 n^{r-2}} \leq \frac{2(r-1) c_3 n^{r-1}}{c_6 n^{r-2}} = c_5 n.
            \end{equation*}
            This implies that
            \begin{align*}
                |N_j(u)| \geq& |V_j| - |X| - c_5 n \\
                        \geq& \frac{3n}{4t} - \frac{r \epsilon_1 n}{c_3} - c_5n \geq c_4 n.
            \end{align*}
        \end{proof}
        
        Let $v_1 v_2$ be a color-critical edge of $F$, and $U_1, U_2, \dots, U_{t}$ be a partition of $F - v_1 v_2$.
        Without loss of generality, we may assume that $v_1, v_2 \in U_1$.
        
        \begin{clm}\label{clm: every vertex has a poor class}
            For every $u \in V(H')$, there exists a $V_i$ for some $i \le t$ that is poor for $u$.
        \end{clm}
        \begin{proof}[Proof of Claim \ref{clm: every vertex has a poor class}]
            Suppose that all the classes are rich for $u$.
            Without loss of generality, let $u \in V_1$.
            We embed $v_1$ to $u$, $v_2$ to $w \in N_1(u) \cap (V_1 \setminus X)$ (this is doable by the $c_3$ and $c_4$ we choose) and $U_1 \setminus \{v_1, v_2\}$ to $V_1 \setminus (X \cup \{u, w\})$.
            
            Next we embed vertices in $U_2$ to $V_2 \setminus X$.
            For each $x \in U_2$, if $x$ is adjacent to $v_1$ in $F$, then the candidate set we choose $x$ from is $N_2(u)$ and $|N_2(u)| \geq c_4 n$.
            If $x$ is nonadjacent to $v_1$ in $F$, then the candidate set we choose $x$ from is $V_2 \setminus X$ and $|V_2 \setminus
            X| \geq \frac{3n}{4t} - \frac{r \epsilon_1 n}{c_3} \geq c_4 n$.
            We also need that the image of $x$ forms heavy pairs with the images of vertices in $U_1$.
            This is doable.
            Since the images of $v_2$ and $U_1$ lie in $V_1 \setminus X$, by Claim \ref{clm: vertex outside X has many heavy partner}, there are at most $|U_1| c_5 n$ vertices that are excluded.
            Thus $c_4 n - |U_1| c_5 n > |U_2|$ for large $n$ and we can embed $U_2$ to $V_2 \setminus X$.
            
            Repeating this process, we embed all vertices of $F$ into $V(H')$.
            Since pairs of vertices in different classes and $u w$ are heavy, we can find an $r$-expansion of $F$ in $H_{\ell}'$ greedily, which corresponds to a rainbow copy of $F^{r+}$, a contradiction.
            This completes the proof.
        \end{proof}
        
        \begin{clm}\label{clm: expansion of color-critical edge}
            Let $i \le t$ and $u \in V(H')$.
            If all the classes except for $V_i$ are rich for $u$, then there is no hyperedge in $H'$ containing $u$ and $v \in V_i \setminus X$.
        \end{clm}
        \begin{proof}[Proof of Claim \ref{clm: expansion of color-critical edge}]
            Suppose otherwise.
            There exists a hyperedge $h$ that contains $u$ and $v$, where $v \in V_i \setminus X$.
            Without loss of generality, assume that $u \in V_j$.
            We embed $v_1$ to $u$, $v_2$ to $v$ and $U_1 \setminus \{v_1, v_2\}$ to $V_i \setminus (X \cup h)$.
            
            Case 1. When $i = j$.
            Similar to the embedding process in the proof of Claim \ref{clm: every vertex has a poor class}, we can embed all the vertices of $F$ into $V(H')$ such that pairs of vertices in different classes are heavy.
            The difference is that we need to avoid vertices in $h$ when we embed vertices in $U_k$, $2\le k \le t$.
            
            Case 2. When $i \neq j$.
            Without loss of generality, assume that $i = 1$ and $j = 2$.
            Next we embed $U_2$ to $V_2 \setminus (X \cup \{u\} \cup h)$.
            For each $x \in U_2$, if $x$ is adjacent to $u$, then the candidate set we choose $x$ from is $N_2(u) \setminus (\{u\} \cup h)$.
            Otherwise, the candidate set we choose $x$ from is $V_2 \setminus (X \cup \{u\} \cup h)$.
            In either case, the size of candidate set is at least $c_4 n - r$.
            We also require that the image of $x$ forms heavy pairs with the images of vertices of $U_1$.
            Similar to the the proof of Claim \ref{clm: every vertex has a poor class}, this is doable.
            For $U_3, \dots, U_{t}$, we can embed them similar to the proof of Claim \ref{clm: every vertex has a poor class}.
            
            In either cases, we can find an $r$-expansion of $F$ such that all hyperedges but $h$ are in $H_{\ell}'$ which corresponds to a rainbow copy of $F^{r+}$, a contradiction.
            This completes the proof.
        \end{proof}
        
        \begin{clm}\label{clm: own class is unique poor class}
            For every $u \in V_i$, the unique poor class for $u$ is $V_i$.
        \end{clm}
        \begin{proof}[Proof of Claim \ref{clm: own class is unique poor class}]
            The statement holds for $u \in V_i \setminus X$ by Claim \ref{clm: vertex outside X has many heavy partner} and Claim \ref{clm: every vertex has a poor class}.
            We assume that $u \in X$.
            First we prove that for every $u \in X$, there is at most one poor class.
            Suppose that $V_p$ and $V_q$ are both poor classes for $u$.
            We count the number of hyperedges in $H_{\ell}'$ containing $u$.
            
            Type 1. Hyperedges containing $u$ and a vertex $v$ where $v \in X$.
            The number of such hyperedges is at most $|X| n^{r-2} \leq \frac{r \epsilon_1}{c_3} n^{r-1}$.
            
            Type 2. Hyperedges containing $u$ and a vertex $v \notin X$ where $(u, v)$ is not a heavy pair.
            The number of such hyperedges is at most $n (r-2)((r-2) \ell + v_F)n^{r-3}$.
            
            Type 3. Hyperedges containing $u$ and a vertex $v \in (V_p \cup V_q) \setminus X$ forming a heavy pair.
            Since $V_p$ and $V_q$ are poor for $u$, the number of such hyperedges is at most $2 c_4 n n^{r-2}$.
            
            Type 4. Hyperedges containing $u$ and all the other $r-1$ vertices are in $\cup_{j \in [t] \setminus \{p, q\}} (V_j \setminus X)$ and form heavy pairs with $u$.
            Since the unique poor class for a vertex not in $X$ is its own class, by Claim \ref{clm: expansion of color-critical edge}, any hyperedge in $H'$ cannot contain two vertices in $V_i \setminus X$ for each $i \le t$.
            Thus, such a hyperedge contains at most one vertex in $V_i \setminus X$ for each $i \le t$.
            So the number of such hyperedges is at most $\binom{t-2}{r-1} \left(\left(\frac{1}{t} + c_2\right) n\right)^{r-1}$.
            
            Summing these four types of hyperedges, we have 
            \begin{align*}
                d_{H'_{\ell}} (u) \leq& \frac{r \epsilon_1}{c_3} n^{r-1} + (r-2)((r-2)\ell + v_F)n^{r-2} + 2 c_4 n^{r-1} + \binom{t-2}{r-1} \left(\left(\frac{1}{t} + c_2\right) n\right)^{r-1} \\
                =& \left(\binom{t-2}{r-1} \left(\frac{1}{t} + c_2\right)^{r-1} + \frac{r \epsilon_1}{c_3} + 2 c_4\right) n^{r-1} + O(n^{r-2}) \\
                \leq& \left(\binom{t-2}{r-1} \left(\frac{1}{t}\right)^{r-1} + \frac{1}{16} \binom{t-2}{r-2} \left(\frac{1}{t}\right)^{r-1} + \frac{r \epsilon_1}{c_3} + 2 c_4\right) n^{r-1} + O(n^{r-2})\\
                =& \left(\binom{t-1}{r-1} \left(\frac{1}{t}\right)^{r-1} - \frac{15}{16} \binom{t-2}{r-2} \left(\frac{1}{t}\right)^{r-1} + \frac{r \epsilon_1}{c_3} + 2 c_4\right) n^{r-1} + O(n^{r-2})\\
                <& |E(T^r(n, t))| - |E(T^r(n-1, t))| - c_1' n^{r-1},
            \end{align*}
            a contradiction.
            So for every $u \in X$, there is exactly one poor class.
            
            Next we prove that the poor class is the class $u$ lies in.
            Suppose otherwise.
            Let $u \in V_i$ and assume that $V_j$ is poor for $u$, $i \neq j$.
            Then we move $u$ from $V_i$ to $V_j$ and denote the new partition by $P'$.
            We need to calculate
            \begin{align*}
                &|E(H_{\ell}' \cap K(P')^r)| - |E(H_{\ell}' \cap K(P)^r)| \\
                =& |E(H_{\ell}' \setminus K(P)^r)| - |E(H_{\ell}' \setminus K(P')^r)| \\
                =& (|\{h \in H_{\ell}' \setminus K(P)^r: u \in h\}| + |\{h \in H_{\ell}' \setminus K(P)^r: u \notin h\}|) \\
                &- (|\{h \in H_{\ell}' \setminus K(P')^r: u \in h\}| + |\{h \in H_{\ell}' \setminus K(P')^r: u \notin h\}|) \\
                =& |\{h \in H_{\ell}' \setminus K(P)^r: u \in h\}| - |\{h \in H_{\ell}' \setminus K(P')^r: u \in h\}| \\
                =& d_{H_{\ell}'}(u) - |\{h \in H_{\ell}' \cap K(P)^r: u \in h\}| - |\{h \in H_{\ell}' \setminus K(P')^r: u \in h\}|.
            \end{align*}
            
            Consider hyperedges in $\{h \in H_{\ell}' \cap K(P)^r: u \in h\}$.
            The number of such hyperedges where $h$ contains a vertex of $X$ other than $u$ is at most $|X| n^{r-2}$.
            If $h$ contains no vertex in $X$, then $h$ does not intersect with $V_j \setminus X$ by Claim \ref{clm: expansion of color-critical edge}.
            So the number of such hyperedges is at most $\binom{t-2}{r-1} ((\frac{1}{t} + c_2)n)^{r-1}$.
            Thus 
            \begin{equation*}
                |\{h \in H_{\ell}' \cap K(P)^r: u \in h\}| \leq \binom{t-2}{r-1} \left(\left(\frac{1}{t} + c_2\right)n\right)^{r-1} + |X| n^{r-2}.
            \end{equation*}
            
            Consider hyperedges in $\{h \in H_{\ell}' \setminus K(P')^r: u \in h\}$.
            Such hyperedges either contain at least two vertices in the same class in partition $P'$ or contain a vertex in $V_j$ in partition $P$.
            For the first case, by the fact that the unique poor class for a vertex outside $X$ is its own class and Claim \ref{clm: expansion of color-critical edge}, these two vertices cannot both be outside $X$.
            This means at least one vertex is in $X$.
            For the second case, by Claim \ref{clm: expansion of color-critical edge}, this vertex is in $X$.
            In either case, each such hyperedge contains a vertex in $X$.
            Thus, the number of such hyperedges is at most $|X| n^{r-2}$.
            Then we have
            \begin{align*}
                &|E(H_{\ell}' \cap K(P')^r)| - |E(H_{\ell}' \cap K(P)^r)| \\
                \geq& d_{{H}_{\ell}'}(u) - \binom{t-2}{r-1} \left(\left(\frac{1}{t} + c_2\right)n\right)^{r-1} - |X| n^{r-2} - |X| n^{r-2} > 0,
            \end{align*}
            a contradiction.
        \end{proof}
        
        \begin{clm}\label{clm: each hyperedge contains two vertices in X}
            Every hyperedge in $H' \setminus K(P)^r$ contains at least two vertices in $X$.
        \end{clm}
        \begin{proof}[Proof of Claim \ref{clm: each hyperedge contains two vertices in X}]
            Let $h \in H' \setminus K(P)^r$.
            Then $h$ contains at least two vertices in some class.
            These two vertices cannot both be outside $X$ by Claim \ref{clm: expansion of color-critical edge} and Claim \ref{clm: own class is unique poor class}.
            It is also impossible that one of these two vertices is in $X$ and the other is outside $X$ by Claim \ref{clm: expansion of color-critical edge}.
            Thus, $h$ must contain at least two vertices of $X$.
        \end{proof}
        Now we can complete the proof.
        \begin{align*}
            e({H}') \leq& k |E(H_{\ell}')| + (\ell - 1) |E(H_1' \setminus H_{\ell}')| \\
            =& k \left(|E(H_{\ell}')| + \frac{\ell - 1}{k} |E(H_1' \setminus H_{\ell}')|\right) \\
            =& k (e(K(P)^r) - |E(H_{\ell}'')| + |E(H_{\ell}' \setminus K(P)^r)| \\
            &+ \frac{\ell - 1}{k} |E((H_1' \setminus H_{\ell}') \cap K(P)^r)| + \frac{\ell - 1}{k} |E((H_1' \setminus H_{\ell}') \setminus K(P)^r))| \\
            \leq& k\left(e(K(P)^r) - \left(1- \frac{\ell - 1}{k}\right)|E(H_{\ell}'')| + |E(H_1' \setminus K(P)^r)|\right).
        \end{align*}
        The last inequality follows from $(H_1' \setminus H_{\ell}') \cap K(P)^r \subseteq H_{\ell}''$.
        For $H_1' \setminus K(P)^r$, by Claim \ref{clm: each hyperedge contains two vertices in X}, we have
        \begin{equation*}
            |E(H_1' \setminus K(P)^r)| \leq \binom{|X|}{2} \binom{n-2}{r-2}.
        \end{equation*}
        For $H_{\ell}''$,
        we have $r |E(H_{\ell}'')| = \sum_{u \in V} d_{H_{\ell}''}(u) \geq \sum_{u \in X} d_{H_{\ell}''}(u) \geq |X| c_3 n^{r-1}$.
        Thus 
        \begin{align*}
            e(H') \leq& k \left(|E(K(P)^r)| - \left(1- \frac{\ell - 1}{k}\right)|E(H_{\ell}'')| + |E(H_1' \setminus K(P)^r)|\right) \\
            \leq& k \left(|E(K(P)^r)| - \left(1- \frac{\ell - 1}{k}\right) |X|\frac{c_3 n^{r-1}}{r} + \binom{|X|}{2} \binom{n-2}{r-2}\right) \\
            \leq& k |E(K(P)^r)| \leq k |E(T^r(n, t))|,
        \end{align*}
        i.e., $f(H') \leq 0$, a contradiction.
        The proof is completed.
    \end{proof}

\section{Concluding remarks}\label{sec4}

AI statement: An incomplete proof of Theorem \ref{thm: sum of edges r-partite} for linear hypergraphs was completed by GPT-6. This was generalized to Theorem \ref{thm: sum of edges r-partite} without use of AI. The authors used GPT-6 to check the calculation in the proof of Theorem \ref{thm: sum of edges expansion of color-critical graph}.

\bigskip
\textbf{Funding}:
The research of Cheng is supported by the National Natural Science Foundation of China (Nos. 12131013 and 12471334), and the China Scholarship Council (No. 202406290241).

The research of Gerbner is supported by the J\'anos Bolyai research scholarship of the Hungarian Academy of Sciences.

The research of Zhang is supported by the China Scholarship Council (No. 202608310129).

\end{document}